\documentclass[11pt]{article}
\usepackage{amssymb}
\usepackage{multicol}
\usepackage{amsmath}
\usepackage{longtable}

\newtheorem{theorem}{Theorem}[section]

\newtheorem{example}{Example}[section]

\begin{document}
\today
\begin{center} 
\begin{Large}
{\LARGE\bf % 138cor/138cora.tex \\
Fredholm equations for non-symmetric kernels,
with applications to  iterated integral operators.
}\\[1ex]
\end{Large}
by\\[1ex]
Christopher S. Withers\\
Applied Mathematics Group\\
Industrial Research Limited\\
Lower Hutt, NEW ZEALAND\\[2ex] Saralees Nadarajah
\\ School of Mathematics\\ University of Manchester\\ Manchester M60 1QD, UK
\end{center}
\vspace{1.5cm}
{\bf Abstract:}~~We give 
the Jordan form and the Singular Value Decomposition
for an integral operator ${\cal N}$ with a non-symmetric kernel $N(y,z)$. 
This is used to give solutions of Fredholm equations for non-symmetric kernels,
 and to determine
 the behaviour of ${\cal N}^n$ and $({\cal N}{\cal N^*})^n$  for large $n$.

\section{Introduction and summary}% section 1
\setcounter{equation}{0}
%\addtocounter{section}{1}

Suppose that $\Omega\subset R^p$ and that $\mu$ is a $\sigma$-finite measure
on $\Omega$.
 Consider a $s_1\times s_2$ complex matrix function $N(y,z)$ on $\Omega\times
\Omega.$
Generally we shall assume that $N\in L_2(\mu\times\mu)$ and is non-trivial,
that is,
$$
0< ||{\cal N}||_2^2= \int\int ||N(y,z)||^2d\mu(y)d\mu(z)<\infty.
$$
The integral operator asociated with $(N,\mu)$ is
${\cal N}$ defined
 by
\begin{eqnarray}
{\cal N}q(y)=\int N(y,z) q(z) d\mu(z),\
p(z)^* {\cal N}
=\int  p(y)^* N(y,z) d\mu(y)
 \label{joint1} 
\end{eqnarray}
where $p:\Omega\rightarrow C^{s_1}$
and  $q:\Omega\rightarrow C^{s_2}$
 are any functions for which these integrals exist, for example
 $p,q\in L_2(\mu)$, that is $\int |p|^2d\mu<\infty$, and similarly for $q$.
(All integrals are over $\Omega.$ * denotes the transpose of the complex conjugate.) 

Section 2 reviews Fredholm theory for
Hermitian kernels, that is, when $N(y,z)^*=N(z,y)$, so that $s_1=s_2$  .
For this case,   ${\cal N}^n=O(r_1^n)$ for large $n$ where $r_1$ is the
magnitude of the largest eigenvalue.

Section 3 extends this to non-Hermitian kernels for the case of diagonal
Jordan form. Again ${\cal N}^n=O(r_1^n)$ for large $n$ for $r_1$ as before.

 Section 4 deals with the  case of non-diagonal
Jordan form. In this case  ${\cal N}^n=O(r_1^n n^{ M-1})$ for large $n$ for $r_1$ as before where $M$ is the largest multiplicity of those eigenvalues
with modulus $r_1$.

Section 5 gives the  Singular Value Decomposition (SVD)  
for a non-symmetric kernel $N(y,z)$. In this case one has results such as
  $({\cal N}{\cal N}^*)^n=O(\theta_1^{2n})
$ for large $n$ where $\theta_1$ is the largest singular value..

\section{Hermitian kernels}% section 2
\setcounter{equation}{0}
%\addtocounter{section}{1}
 {\bf Matrix theory}\\
First consider a Hermitian matrix $N^*=N\in C^{s\times s}.$
Its eigenvalues $\nu_1,\cdots,\nu_s$ are the roots of $\det(N-\nu I)=0.$
They are real.
Corresponding to $\nu_j$ is an eigenvector $p_j$ satisfying $Np_j=\nu_jp_j$.
These are orthonormal: $p_j^*p_k=\delta_{jk}$ where
$\delta_{jj}=1$ and $\delta_{jk}=0$ for $j\neq k$. 
Set $P=(p_1,\cdots,p_s)$. If $N$  and its eigenvalues are real, then  $P$ can be taken
to be real.
The {\it spectral decomposition} of $N$ in terms of its eigenvalues and eigenvectors is
\begin{eqnarray}
N &=& P\Lambda P^*=\sum_{j=1}^s \nu_j p_jp_j^* \mbox{ where }\Lambda=\mbox{diag}(\nu_1,\cdots, \nu_s), \label{herm} \\
\sum_{j=1}^s \nu_j p_jp_j^* &=&P P^*= I_s=P^*P=(p_j^*p_k).\nonumber
\end{eqnarray}
So for $\alpha\in C$
\begin{eqnarray}
N^\alpha=P\Lambda^\alpha P^*=\sum_{j=1}^s \nu_j^\alpha p_jp_j^* 
, \label{alpha} 
\end{eqnarray}
provided that if $\det(N)= 0$, then $\alpha$ has non-negative real part.
So for $\det(N)\neq 0$, 
$$Nf=g\Rightarrow f=N^{-1}g=\sum_{j=1}^s \nu_j^{-1}p_j(p_j^*g).$$
Similarly for $\nu$ not an eigenvalue and $f,g\in C^s$,
$$(\nu I-N)f=g\Rightarrow f=(\nu I-N)^{-1}g=\sum_{j=1}^s (\nu-\nu_j)^{-1}p_j(p_j^*g).$$
For large $n$, if 
\begin{eqnarray}
r_1=|\nu_1|=\cdots=|\nu_M|>r_0=\max_{j=M+1}^s |\nu_j|
 \label{M} 
\end{eqnarray}
 then
\begin{eqnarray}
N^n=r_1^nC_n+O(r_0^n) 
\mbox{ where }C_n=\sum_{j=1}^M [\mbox{sign }(\nu_j)]^n p_jp_j^* = O(1)
 \label{hermapprox} 
\end{eqnarray}
assuming that $s$ does not depend on $n$.

 {\bf Function theory}\\
Now consider  a function $N(y,z):\Omega^2\rightarrow C^{s\times s}$
and a $\sigma$-finite measure $\mu$ on $\Omega\subset R^p$.
Its integral operator with respect to $\mu$ is ${\cal N}$ defined by 
(\ref{joint1}). Suppose that the {\it kernel} $N$ is Hermitian, that is, 
$$N(y,z)^*=N(z,y).$$ Then the analogues of the matrix results above are
as follows. 
Suppose that $ ||{\cal N}||_2^2>0$ and 
$$\sum_{j=1}^s\int |N_{jj}(x,x)d\mu(x)<\infty.$$
The {\it spectral decomposition} of $N$ in terms of its eigenvalues and 
vector eigenfunctions $\{p_j(y)\}:\Omega\rightarrow C^s$, is
\begin{eqnarray}
N(y,z) &=& P(y)\Lambda P(z)^*=\sum_{j=1}^\infty \nu_j p_j(y)p_j(z)^* \label{dvj} \\
\mbox{where }\Lambda&=& \mbox{diag}(\nu_1,\nu_2,\cdots),\
P(y)=(p_1(y),p_2(y),\cdots)
. \nonumber
\end{eqnarray}
The eigenfunctions are orthonormal with respect to $\mu$:
$$\int p_j^*p_kd\mu=\delta_{jk}.$$
In Fredholm theory the convention is to call $\{\lambda_j=\nu_j^{-1}\}$
the eigenvalues, rather than $\{\nu_j\}$.

${\cal I}_{yz}=P(y) P(z)^*$ is generally divergent but can be thought of as a
generalized Dirac function: $\int {\cal I}_{yz} f(z)d\mu(z)=f(y).$
 $P(z)^*P(y)=(p_j(z)^*p_k(y))$ satisfies  $\int P(z)^*P(z)d\mu(z)=I_\infty.$

 For $n=1,2,\cdots, \ N_n(y,z)={\cal N}^{n-1}N(y,z)$ satisfies
\begin{eqnarray*}
 N_n(y,z)=P(y)\Lambda^n P(z)^*=\sum_{j=1}^s \nu_j^n p_j(y)p_j(z)^*.
\end{eqnarray*}
For large $n$, if (\ref{M}) holds then
\begin{eqnarray*}
N_n(y,z)=r_1^nC_n(y,z)+O(r_0^n) 
\mbox{ where }C_n(y,z)=\sum_{j=1}^M [\mbox{sign }(\nu_j)]^n p_j(y)p_j(z)^*=O(1).
 \label{hermapproxfn} 
\end{eqnarray*}
 Conditions for (\ref{dvj}) to hold pointwise
and uniformly are given in Withers (1974, 1975, 1978.)
It is known as Mercer's Theorem.

{\bf The resolvent}\\
Given functions $f,g:\Omega\rightarrow C^s$, 
the {\it Fredholm integral equation of the second kind},
\begin{eqnarray}
p(y)-\lambda{\cal N}p(y)=f(y),\label{r1}
\end{eqnarray}
can be solved for $\lambda$ not an eigenvalue using
$$ 
(I-\lambda{\cal N})^{-1}=I+\lambda{\cal N}_\lambda,
\mbox{ that is, }
{\cal N}_\lambda=(I-\lambda{\cal N})^{-1}{\cal N}
$$
where
$${\cal N}_\lambda f(y)=\int N_\lambda(y,z)f(z)d\mu(z),\
g(z){\cal N}_\lambda=\int g(y)N_\lambda(y,z)d\mu(y),
$$
and the {\it resolvent}  of ${\cal N}$,
$$N_\lambda(y,z)=(I-\lambda{\cal N})^{-1} N(y,z)
:{\cal C}\times\Omega^2
\rightarrow {\cal C}^{s\times s}
$$ 
with operator ${\cal N}_\lambda$
is the unique solution of
$$ (I-\lambda{\cal N}){\cal N}_\lambda={\cal N}
={\cal N}_\lambda(I-\lambda{\cal N}),
$$
that is,
$$\lambda {\cal N} N_\lambda(y,z)= N(y,z)- N_\lambda(y,z)
=\lambda  N_\lambda(y,z){\cal N}.
$$
If this can be solved analytically or numerically, then one has a solution of
(\ref{r1}) without the
need to compute the eigenvalues and eigenfunctions of ${\cal N}$.

The resolvent satisfies
\begin{eqnarray}
N_\lambda(y,z)=\sum_{j=1}^\infty p_j(y)p_j(z)^*/(\lambda_j-\lambda).
\label{res2}
\end{eqnarray}
 Conditions for this to hold are given by Corollary 3 of Withers (1975).
 The Fredholm equation of the second kind, (\ref{r1}), has solution
$$
p(y)=f(y)+\lambda\sum_{j=1}^\infty p_j(y)\int p_j^*fd\mu/(\lambda_j-\lambda).
$$
The resolvent exists except for $\lambda=\lambda_j$, an eigenvalue.
The eigenvalues of ${\cal N}$ are the zeros of its
 {\it Fredholm determinant}
\begin{eqnarray}
D(\lambda)=\Pi_{j=1}^\infty (1-\lambda/\lambda_j)
=\exp\{-\int_0^\lambda d\lambda\int \mbox{trace }N_\lambda(x,x)d\mu(x)\}.
 \label{fred} 
\end{eqnarray}
The {\it Fredholm integral equation of the first kind}
$$ \lambda {\cal N}p(x)= p(x)$$
has a solution provided that $\lambda$ is an eigenvalue. 
For $\nu$ an eigenvalue, its general solution $p(x)$ is a linear combination
of the eigenfunctions $\{p_j(x)\}$ corresponding to $\lambda_j=\lambda$. 
\begin{example} %2.1
Suppose that $Y,Z\in R$ and ${Y\choose Z}\sim{\cal N}_{2}(0,V)$
where $V= \begin{pmatrix} I & r\\ r & I\end{pmatrix}$.
So $V$ is the correlation matrix for ${Y\choose Z}$.
For $j\in N_+,\ x\in R$ set $p_j(x)=H_{j}(x)/j!^{1/2}$ where
 $ H_j(x)$ is the standard univariate Hermite polynomial. Then
$\int p_jp_k\phi_I=\delta_{jk}$ and
$$\sum_{j=0}^\infty r^{j}p_j(y)p_j(z)=\phi_C(y,z)/\phi_I(y)\phi_I(z).$$
This is  Mehler's expansion for the standard bivariate normal distribution.
Pearson gave an integrated version and 
Kibble extended it to an expansion for
$\phi_V(x)/\phi_I(x)$ for $x\in R^k$ and $V$ a correlation matrix. 
See for 
example (45.52) p127 and p321-2 of Kotz, Balakrishnan and Johnson (2000). 
\end{example} %2.1
%\begin{example} %2.2 Give resultant of W and $W_0$K for prev. section? \end{example} %2.2

\section{Functions of two variables with diagonal Jordan form}% section 3
\setcounter{equation}{0}
%\addtocounter{section}{1}
{\bf Diagonal Jordan form for matrices}\\
Consider $N\in C^{s\times s}$ with eigenvalues $\nu_1,\cdots, \nu_s$,
the roots of $\det(N-\nu I)=0.$ $N$ is said to have {\it diagonal Jordan form}
 (DJF) if
$$N=P\Lambda Q^* \mbox{ where }PQ^*=I
\mbox{ and }
\Lambda=\mbox{diag}(\nu_1,\cdots, \nu_s).
$$ % as in {herm} 
So 
\begin{eqnarray*}
N &=& \sum_{j=1}^s \nu_jp_jq_j^*,\ q_j^*p_k=\delta_{jk}\mbox{ where }
P=(p_1,\cdots, p_s),\ Q=(q_1,\cdots, q_s),\\
\sum_{j=1}^s p_jq_j^* &=& P Q^*=I_s= Q^*P=(q_j^*p_k).
\end{eqnarray*}
If $N,\Lambda$ are real then $P,Q$ can be taken as real.
Also
\begin{eqnarray}
NP=P\Lambda,\ Np_j=\nu_jp_j,\
N^*Q=Q{\overline \Lambda},\ N^*q_j={\overline \nu}_j q_j,\label{cfp10}
\end{eqnarray}
and for any complex $\alpha$,
$$N^\alpha=P\Lambda^\alpha Q^*=\sum_{j=1}^s \nu_j^\alpha p_jq_j^*
$$
provided that if $\det(N)= 0$, then $\alpha$ has non-negative real part.
Suppose that
\begin{eqnarray}
\nu_j=r_je^{i\theta_j}\mbox{ and }r_1=\cdots=r_R>r_0=\max_{j=R+1}^s r_j.
 \label{r} 
\end{eqnarray}
Then
\begin{eqnarray}
N^n=r_1^nC_n+O(r_0^n) 
\mbox{ where }C_n=\sum_{j=1}^R e^{in\theta_j} p_jq_j^* = O(1)
 \label{hermapprox2} 
\end{eqnarray}
Taking $\alpha=-1$ gives the  inverse of $N$ when this exists. 

 {\bf Diagonal Jordan form for functions}\\
Now consider  a function $N(y,z):\Omega^2\rightarrow C^{s\times s}$.
When $N$ has diagonal Jordan form (for example when its eigenvalues are all
different), then
{\it the Fredholm equations of the first kind},
$$
\lambda {\cal N}p(y)=p(y),\ 
\bar{\lambda}{\cal N}^* q(z)= q(z), 
$$
or  equivalently for $\nu=\lambda^{-1}$,
$$
{\cal N}p(y)=\nu p(y),\ 
{\cal N}^* q(z)= \bar{\nu} q(z), 
$$
also have only a countable number of solutions, say
  $\{\lambda_{j}=\nu_j^{-1},p_{j}(y),q_{j}(z),\ j\geq 1 \}$
up to {\it arbitrary} constant multipliers for   $\{p_{j}(y),\ j\ge 1 \}$,
 satisfying the {\it bi-orthogonal} conditions
$$\int q_{j}^*p_{k}d\mu=\delta_{jk}.
$$
These are called the {\it eigenvalues} and {\it right and left eigenfunctions}
of $(N,\mu)$ or  ${\cal N}$.
Also
\begin{eqnarray}
 N(y,z) &=& \sum_{j=1}^\infty \nu_j p_j(y)q_j(z)^*
=P(y)\Lambda Q(z)^*
 \label{djf} \\
\mbox{where }\Lambda&=& \mbox{diag}(\nu_1,\nu_2,\cdots),\
P(y)=(p_1(y),p_2(y),\cdots),\ Q(z)=(q_1(y),q_2(y),\cdots)
. \nonumber
\end{eqnarray}
with convergence in $L_2(\mu\times\mu)$, or pointwise and uniform under 
stronger conditions.
If $N$ is a real function and $\Lambda$ is real,
then $P,Q$ can be taken as real functions

 For $n\geq 1$,
\begin{eqnarray}
 N_n(y,z)={\cal N}^{n-1}N(y,z)
 \label{Nn} 
\end{eqnarray}
satisfies
\begin{eqnarray}
 N_n(y,z) &=& \sum_{j=1}^\infty \nu_j^n p_j(y)q_j(z)^*=P(y)\Lambda^n Q(z)^*,
\label{Nnp} \\
 {\cal N}^{n}p(y) &=& \sum_{j=1}^\infty \nu_j^n p_j(y)\int q_j^*pd\mu. 
\nonumber
\end{eqnarray}
If (\ref{r}) holds then
\begin{eqnarray}
N_n(y,z)&=& r_1^nC_n+O(r_0^n) 
\mbox{ where }C_n=\sum_{j=1}^R e^{in\theta_j} p_j(y)q_j(z)^* = O(1)
 \label{hermapproxfn2}  \\
 {\cal N}^{n}p(y) &=& r_1^nc_n+O(r_0^n) 
\mbox{ where }c_n=\sum_{j=1}^R  e^{in\theta_j} p_j(y)\int q_j^*pd\mu = O(1). 
\nonumber
\end{eqnarray}
The resolvent satisfies the equations of Section 2 except that 
(\ref{res2}) is replaced by
\begin{eqnarray}
N_\lambda(y,z)=\sum_{j=1}^\infty p_j(y)q_j(z)^*/(\lambda_j-\lambda).
\label{res3}
\end{eqnarray}
The {\it Fredholm determinant} is again given by (\ref{fred}).
If only a finite number of eigenvalues are non-zero, the {\it kernel} $N(y,z)$
is said to be {\it degenerate}.
(For example this holds if $\mu$ puts weight only at $n$ points.)
 If $R=1$, that is,
\begin{eqnarray}
|\lambda_1|<|\lambda_j|\mbox{ for }j>1, \label{l1}
\end{eqnarray}
then as $n\rightarrow\infty$,
$$
{\cal N}^{n+1} f(y)/{\cal N}^n f(y)\rightarrow\lambda_1^{-1},\
 f(y){\cal N}^{n+1}/ f(y){\cal N}^n \rightarrow\lambda_1^{-1}.
$$
This is one way to obtain the first eigenvalue
$\lambda_1$ arbitrarily closely. Another is to use
\begin{eqnarray}
\lambda_{1}^{-1}=\sup\{ \int g {\cal N} h d\mu:\ \int ghd\mu=1\}
\mbox{ if }\lambda_{1}>0,
 \label{1}\\
\lambda_{1}^{-1}=\inf\{ \int g {\cal N} h d\mu:\ \int ghd\mu=1\}
\mbox{ if }\lambda_{1}<0.
 \label{2}
\end{eqnarray}
The maximising/minimising functions are the first eigenfunctions
$g=g_1,h=h_1$. These are unique up to a constant
multiplier if (\ref{l1}) holds.
If  $\lambda_1$ is known, one can use
$$
(\lambda_1{\cal N})^{n} f(y)\rightarrow p_1(y)\int q_1^*fd\mu,\
 f(y)^*(\lambda_1{\cal N})^{n}\rightarrow q_1(y)\int f^*p_1d\mu,\
$$
for any function $f:\Omega\rightarrow {\cal C}^s$,
to approximate $p_1(y), q_1(y)$.
%Also since $q_1(y)$ is only unique up to a multiplicative constant, we may choose say $\int q_1^* f=1$ and so approximate  $p_1(y), q_1(y)$.
One may now repeat the procedure on the operator
${\cal N}_1$ corresponding to
$$N_1(y,z)=N(y,z)-\nu_1 p_1(y)q_1(z)^*$$
to approximate $\lambda_{2},p_{2}(y),q_{2}(z)$
assuming that the next eigenvalue in magnitude, $\lambda_{2},$ has multiplicity 1.
%If %say  $\lambda_{1}$  has multiplicity $R>1$, then $$ (\lambda_1{\cal N})^{n} f(y)\rightarrow \sum_{j=1}^R p_j(y)\int q_j^*fd\mu,$$ and one can adapt the method above.

For further details see Withers (1974, 1975, 1978) and references.

\section{Fredholm theory for non-diagonal Jordan form}% section 4
\setcounter{equation}{0}
%\addtocounter{section}{1}
 
{\bf Non-diagonal Jordan form for matrices}

For $N\neq N^*$ a matrix in $C^{s\times s}$, its general Jordan form is
\begin{eqnarray}
N &=& PJP^{-1}
\mbox{ where }J=diag(J_1,\cdots,J_r),\ 
J_j=J_{m_j}(\lambda_j),\ \sum_{j=1}^r m_j=s,\nonumber \\
 J_{m}(\lambda) &=& \lambda I_m+ U_m
=
 \begin{pmatrix}
\lambda & 1 & 0 & \cdots & 0 \\
0 &  \lambda     & 1 & \cdots & 0 \\
  &        &        & \cdots &           \\
0 &  0     &  0     &\cdots  &  \lambda
\end{pmatrix}
,\label{gen} 
\end{eqnarray}
for some matrix $P$, and $U_m$ is the $m\times m$ matrix with 1s on the
superdiagonal and 0s elsewhere:
$$(U_m)_{jk}=\delta_{j,k-1}.$$
(See [1] for example.
 If $N$ and its eigenvalues are real, then $P$ can be taken as real.)
So for $n\geq 1$,
$$N^n=PJ^nP^{-1}\mbox{ where }J^n=diag(J_1^n,\cdots,J_r^n).$$
 By the  Binomial Theorem,
$$ J_{m}(\lambda)^n=\sum_{a=0}^n {n\choose a}\lambda^{n-a}U_m^a
\mbox{ and } 
(U_m^a)_{jk}=\delta_{j,k-a}.$$
So $U_m^m=0$.
For example
\begin{eqnarray*}
 J_{2}(\lambda)^n=\lambda^{n}I_2 +n\lambda^{n-1}U_2
= \begin{pmatrix}
\lambda^n & n\lambda^{n-1}\\
0 &  \lambda^n
\end{pmatrix}
.
\end{eqnarray*}
So $N^n$ can be expanded in block matrix form
$$
 (N^n)_{jk} = 
\sum_{c=1}^r P_{jc}\ J_c^n\  P^{ck} 
$$
where we partition $P$ and its inverse as
$$P=(P_{jk}:\ j,k=1,\cdots,r),\ P^{-1}=(P^{jk}:\ j,k=1,\cdots,r)$$
 with elements $P_{jk}$ and  $P^{jk}$
matrices in $C^{m_j\times m_k}$.

Alternatively setting 
$$Q^*=P^{-1},\ (P_1,\cdots, P_r)=P,\  (Q_1,\cdots, Q_r)=Q,$$
with $P_j,Q_j\in C^{s\times m_j}$, we have
\begin{eqnarray}
N^n &=& PJ^nQ^*=\sum_{j=1}^r P_jJ_j^nQ_j^*,\label{alt} \\
\sum_{j=1}^r P_jQ_j^* &=& PQ^*=I_s=Q^*P=(Q_j^*P_k)\nonumber
\end{eqnarray}
so that
$$Q_j^*P_j=I_{m_j},\
Q_j^*P_k=0\in C^{m_j\times m_k}\mbox{ if }j\neq k.
$$
$P$ can be obtained as follows. Let $p_{jk}$ be the $k$th column of $P_j$ for
$k=1,\cdots,m_j$.
 Then
\begin{eqnarray}
NP=PJ\Rightarrow NP_j=P_jJ_j\Rightarrow Np_{jk}=\lambda_jp_{jk}+p_{j,k-1}
\mbox{ where }p_{j0}=0.
\label{alt1} 
\end{eqnarray}
 So one first obtains $p_{j1}$, the right eigenvector of $N$,
then $p_{j2},\cdots,p_{jm_j}$. This is called {\it the Jordan chain}.
$Q$ can either be obtained by inverting $P$ or using
\begin{eqnarray}
N^*Q=QJ^*\Rightarrow N^*Q_j=Q_jJ_j^*\Rightarrow
 N^*q_{jk}=\bar{\lambda}_jq_{jk}+q_{j,k+1}\mbox{ where }q_{j,m_j+1}=0.
\label{alt2} 
\end{eqnarray}
So one first computes $q_{j,m_j}$, the right eigenvector of $N^*$ then
 $q_{j,m_j-1},\cdots,q_{j1}$.
For large $n$ and $\lambda\neq 0$, 
$$J_n(\lambda)^n= {n\choose m-1}\lambda^{n-m+1}[U_m^{m-1}+O(1)]$$
and $U_m^{m-1}$ is a matrix of 0's except for a 1 in its upper right corner.
So if (\ref{r}) holds and
$$m_1=\cdots =m_M>\max_{j=M+1}^R m_j,$$
then %for $1\leq j\leq M$, $J_j^n\approx  {n\choose M-1}\lambda^{n-M+1}U_M^{M-1}$$
\begin{eqnarray*}
 (N^n)_{jk} &=&  {n\choose M-1}r_1^{n-M+1}[D_n+O(n^{-1})]\\
\mbox{where }
D_n &=& \sum_{c=1}^M P_{jc}\ e^{i(n-M+1)\theta_c} U_M^{M-1} P^{ck} 
= \sum_{c=1}^M P_{jc}\ e^{i(n-M+1)\theta_c} U_M^{M-1} P^{ck} =O(1).
\end{eqnarray*}
See Withers and Nadarajah (2008) for more details.

 {\bf Non-diagonal Jordan form for functions}\\
Now consider $N:\Omega^2\rightarrow C^{s\times s}$.
Suppose that $\mu$ is a $\sigma$-finite measure
on $\Omega$ and that $N$ is not Hermitian, that is $N(y,z)^*\neq N(z,y).$
Its Jordan form is
\begin{eqnarray}
N(y,z)=P(y)JP(z)^{-1}=
\mbox{ where }J=diag(J_1,J_2,\cdots),\ 
J_j=J_{m_j}(\lambda_j)
\label{fn} 
\end{eqnarray}
for $P(y):\Omega\rightarrow C^{s\times \infty}$ and
$ J_{m}(\lambda)$ of (\ref{gen}) above. So partitioning 
$$
P(y)=(P_{jk}(y):\ j,k=1,2,\cdots),\ P(z)^{-1}=(P^{jk}(z):\ j,k=1,2,\cdots),
$$
 with elements $P_{jk}(y)$ and  $P^{jk}(z)$
matrix functions in $\Omega\rightarrow C^{m_j\times m_k}$,
we can partition the $n$th iterated kernel, 
$N_n(y,z)={\cal N}^{n-1}N(y,z)$ as
$$[N_n(y,z)]_{jk}=\sum_{c=1}^\infty P_{jc}(y)\ J_c^n\  P^{ck}(z).$$ 
Alternatively setting 
$$Q(z)^*=P(z)^{-1},\ (P_1(y), P_2(y),\cdots)=P(y),\  (Q_1(z), Q_2(z),\cdots)=Q(z),$$
with $P_j(y),Q_j(z):\Omega\rightarrow C^{s\times m_j}$, we have
\begin{eqnarray}
N_n(y,z) &=& P(y)J^nQ(z)^*=\sum_{j=1}^\infty P_j(y)J_j^nQ_j(z)^*.\label{alt3} 
\end{eqnarray}
$P(y)$ can be obtained as follows. Let $p_{jk}(y)$ be the $k$th column of $P_j(y)$ for
$k=1,\cdots,m_j$.
 Then
\begin{eqnarray}
{\cal N}P(y)
=P(y)J\Rightarrow{\cal N}P_j(y)=P_j(y)J_j\Rightarrow{\cal N}p_{jk}(y)=\lambda_jp_{jk}(y)+p_{j,k-1}(y)
\label{alt31} 
\end{eqnarray}
where $p_{j0}(y)=0$. So one first obtains $p_{j1}(y)$, the right eigenfunction of $N$,
then $p_{j2}(y),\cdots,p_{jm_j}(y)$.
$Q(z)
$ can either be obtained by inverting $P(z)$ or using
\begin{eqnarray}
{\cal N}^*Q(z)=Q(z)J^*\Rightarrow {\cal N}^*Q_j(z)=Q_j(z)J_j^*\Rightarrow
 {\cal N}^*q_{jk}(z)=\bar{\lambda}_jq_{jk}(z)+q_{j,k+1}(z)
\label{alt4} 
\end{eqnarray}
 where $q_{j,m_j+1}(z)=0.$
So one first computes $q_{j,m_j}(z)$, the right eigenfunction of $N^*$ then
 $q_{j,m_j-1}(z),\cdots,q_{j1}(z)$.

So if (\ref{r}) holds and
$$m_1=\cdots =m_M>\max_{j=M+1}^R m_j,$$
then %for $1\leq j\leq M$, $J_j^n\approx  {n\choose M-1}\lambda^{n-M+1}U_M^{M-1}$$
\begin{eqnarray*}
 (N_n(y,z))_{jk} &=&  {n\choose M-1}r_1^{n-M+1}(D_{jkn}(y,z)+O(n^{-1}))\\
\mbox{where }
D_{jkn}(y,z) &=& \sum_{c=1}^M \ e^{i(n-M+1)\theta_c}P_{jc}(y) U_M^{M-1} P^{ck}(z) 
 =O(1)
\end{eqnarray*}
has $(a,b)$ element 
$\sum_{c=1}^M \ e^{i(n-M+1)\theta_c}[P_{jc}(y)]_{a1} [P^{ck}(z)]_{Mb}.$ 
\begin{example} %4.1
\end{example} %4.1

\section{The SVD for functions of two variables}% section 4?
\setcounter{equation}{0}
%\addtocounter{section}{1}

{\bf The SVD for matrices}\\
Suppose that $N\in C^{s_1\times s_2}.$ That is, $N$ is
a $s_1\times s_2$ complex matrix. Denote its complex conjugate transpose by $N^*$.
Its SVD is
\begin{eqnarray}
N &=& PD Q^*=\sum_{j=1}^r \theta_j p_j q_j^* \mbox{ where }PP^*=I,\ QQ^*=I,
\ r=\min(s_1,s_2),\label{svdm}\\ 
P &=& (p_1,\cdots,p_{s_1})\in C^{s_1\times s_1},\ 
Q=(q_1,\cdots,q_{s_2})\in C^{s_2\times s_2},
\ \theta_1\geq \cdots\geq \theta_r> 0 \nonumber
\end{eqnarray}
%CONFLICT IN $\theta$ NOTATION FROM DJN Section.\\
and for $s_1=s_2,\ s_1>s_2,\ s_1<s_2$
$$
D=\Lambda,\ {\Lambda\choose 0},\ (\Lambda,0)\mbox{ respectively where }
\Lambda=\mbox{diag}(\theta_1,\cdots,\theta_r).$$
If $N$ is real, then so are $P$ and $Q$.

So for  $s_1>s_2$, 
$$
DD^*={\Lambda^2\ 0\choose 0\ 0},\
D^*D=\Lambda^2
$$
and
 for  $s_1<s_2$, 
$$
DD^*=\Lambda^2,\
D^*D={\Lambda^2\ 0\choose 0\ 0}
.$$
Compare this with (\ref{cfp10}).
Also for $1\leq j\leq r$,
$$N q_j=\theta_j p_j,\ N^* p_j=\theta_j q_j,$$
 for $r< j\leq s_1,\ N q_j=0,$
and for $r< j\leq s_2,\ N^* p_j=0.$
Also since
$$
NN^*P=PDD^*,\ N^*NQ=QD^*D,
$$
the $p_j$ is a right eigenvector of $NN^*$ with eigenvalue 
$\theta_j^2$ (or 0 if  $r< j\leq s_1$) and
the $q_j$ is a right eigenvector of $N^*N$ with 
 eigenvalue $\theta_j^2$ (or 0 if  $r< j\leq s_2$). So (or by Section 2),
\begin{eqnarray}
(NN^*)^n &=& \sum_{j=1}^r \theta_j^{2n} p_jp_j^*,\ 
(N^*N)^n=\sum_{j=1}^r \theta_j^{2n} q_jq_j^*\mbox{ for }n\geq 1,\nonumber\\
(NN^*)^nN &=& \sum_{j=1}^r \theta_j^{2n+1} p_jq_j^*,\ 
(N^*N)^nN^*=\sum_{j=1}^r \theta_j^{2n+1} q_jp_j^*\mbox{ for }n\geq 0.
 \label{nn*} 
\end{eqnarray}
These do not depend on the vectors $\{p_j,q_j,j\geq r\}$. 

So if $\theta_1=\cdots=\theta_M>\theta_{M+1}$, then we have the approximations 
 for $n\geq 0$,
\begin{eqnarray}
(NN^*)^n &=&  \theta_1^{2n}
\sum_{j=1}^M
 p_jp_j^* +O(\theta_{M+1}^{2n}),\ 
(N^*N)^n= \theta_1^{2n}\sum_{j=1}^M q_jq_j^* +O(\theta_{M+1}^{2n})
\mbox{ for }n\geq 1,\nonumber\\
(NN^*)^nN &=& 
 \theta_1^{2n+1}
\sum_{j=1}^r p_jq_j^* +O(\theta_{M+1}^{2n+1}),\ 
(N^*N)^nN^*= \theta_1^{2n+1}
\sum_{j=1}^r q_jp_j^* +O(\theta_{M+1}^{2n+1}).
 \label{apnn*} 
\end{eqnarray}
But 
$$ I_{s_1}=(NN^*)^0 = \sum_{j=1}^{s_1} p_jp_j^*,\ 
I_{s_1}=(N^*N)^0=\sum_{j=1}^{s_2} q_jq_j^*.
$$
If $s_1=s_2$ and $N$ is non-singular, its inverse is
$$N^{-1}=Q\Lambda^{-1} P^*.$$
%If it is singular, a pseudoinverse is given by $$N^{-}=Q\Lambda^{-} P^*.$$
However unlike Jordan form, the SVD does not give a nice form for powers of
$N$.

Now suppose $\Omega\subset R^p$ and
that $\mu$ is a $\sigma$-finite measure on $\Omega$.
 Consider a function $N(y,z):\Omega^2\rightarrow C^{s_1\times s_2}.$
%\\ \ [[DELETE and that $$ ||{\cal N}||_2^2=  \int\int  \sum_{j=1}^{s_1}  \sum_{k=1}^{s_2} N(y,z)N(z,y)d\mu(y)d\mu(z)<\infty,\ \int  \sum_{j=1}^{s} N_{jk}(y,z)N_{kj}(z,y)d\mu(y)d\mu(z)<\infty. FIX!$$ \ ]]\\

 The equations
%$$ \lambda{\cal N}p(y)=p(y),\ \bar{\lambda} q(z){\cal N}^*= q(z), $$ or equivalently putting $\theta=\lambda^{-1}$,
$$
{\cal N}q(y)=\theta p(y),\ 
{\cal N}^* p(z)= {\theta} q(z), 
$$
have a countable number of solutions, say
  $\{\theta_{j},p_{j}(y),q_{j}(z),\ j\geq 1 \}$
 satisfying
$$\int p_{j}^*p_{k}d\mu=\int q_{j}^*q_{k}d\mu=
\delta_{jk}.
$$
The {\it singular values } $\{\theta_{j}\}$ may be taken as real, non-negative
 and
non-increasing. (For convenience we have included $\theta_j=0$.)
$\{p_{j}(y)\}$ and $\{q_{j}(y)\}$
 are the right eigenfunctions of   ${\cal NN^*}$ and  ${\cal N^*N}$ 
respectively, with eigenvalues $\{\theta_j^2\}$.
Also in $L_2(\mu\times\mu)$
\begin{eqnarray}
 N(y,z)=\sum_{j=1}^\infty \theta_j p_j(y)q_j(z)^*. \label{svd} 
\end{eqnarray}
%This is the functional form of the SVD for a square non-symmetric matrix. 
 If $N$ is real, then so are $\{p_j, q_j\}.$
By (\ref{svd} ), for $n\geq 0,$
\begin{eqnarray*}
 ({\cal N}{\cal N}^*)^n N(y,z) &=& \sum_{j=1}^\infty \theta_j^{2n+1} p_j(y)q_j(z)^*, \\
({\cal N}^*{\cal N})^n N(y,z)^* &=& \sum_{j=1}^\infty \theta_j^{2n+1} q_j(z)p_j(y)^*, \\
{\cal N}^* ({\cal N}{\cal N}^*)^n N(y,z) &=& \sum_{j=1}^\infty \theta_j^{2n+2} q_j(y)q_j(z)^*,\\ 
{\cal N}({\cal N}^*{\cal N})^n N(y,z)^* &=& \sum_{j=1}^\infty \theta_j^{2n+2} p_j(z)p_j(y)^*,
\end{eqnarray*}
\begin{eqnarray*}
({\cal N}{\cal N}^*)^n p(y) 
&=&  \sum_{j=1}^\infty \theta_j^{2n}p_j(y)\int p_j^*pd\mu ,\\
({\cal N}^*{\cal N})^n q(z) 
&=&  \sum_{j=1}^\infty \theta_j^{2n}q_j(z)\int q_j^*qd\mu ,\\
({\cal N}{\cal N}^*)^n{\cal N} q(y) 
&=&  \sum_{j=1}^\infty \theta_j^{2n+1}p_j(y)\int q_j^*qd\mu ,\\
({\cal N}^*{\cal N})^n{\cal N}^* p(y) 
&=&  \sum_{j=1}^\infty \theta_j^{2n+1}q_j(y)\int p_j^*pd\mu,\\
\int \mbox{trace } {\cal N}^* ({\cal N}{\cal N}^*)^n N(y,z)|_{z=y}d\mu(y) 
&=&\int \mbox{trace } {\cal N} ({\cal N}^*{\cal N})^n N(y,z)|_{z=y}d\mu(y) 
= \sum_{j=1}^\infty \theta_j^{2n+2} .
% \label{x} 
\end{eqnarray*}
% OR USE $N(z,y)^*,\ p_j(y)p_j(z)^*$ ??
So if $\theta_1=\cdots=\theta_M>\theta_{M+1}$, then we have approximations 
such as
$$
({\cal N}{\cal N}^*)^n p(y) 
=  \theta_1^{2n} \sum_{j=1}^\infty p_j(y)\int p_j^*pd\mu 
+O(\theta_{M+1}^{2n})
$$ if $ p_j(y)\int p_j^*pd\mu =O(1)$ for $j>M$.
So for iterations of ${\cal N}^*{\cal N}$  or  ${\cal N}{\cal N}^*$ 
the most important parameter is the largest singular value.


\begin{thebibliography}{999}

\bibitem{}  Kotz,S., Balakrishnan, N. and Johnson, N.L. (2000)  {\it 
Continuous multivariate distributions}, {\bf 1}, 2nd edition, Wiley, New York.
 
\bibitem{}
Withers, C. S. (1974) Mercer's Theorem and Fredholm resolvents.
         {\it Bull. Austral. Math. Soc.}, {\bf 11},  373-380. %cf ms 94,96

\bibitem{}
 Withers, C. S. (1975) Fredholm theory for arbitrary measure spaces.
 {\it Bull. Austral. Math. Soc.}, {\bf 12}, 283-292.

\bibitem{}
 Withers, C. S. (1978) Fredholm equations have uniformly convergent solutions.
 {\it Jnl. of . Math. Anal. and Applic.}, {\bf 64}, 602-609.
 
\bibitem{}
Withers, C. S. (2000) A simple expression for the multivariate Hermite polynomials.
 {\it Statistics and Probability Letters},  {\bf 47}, 165-169.	%215

\bibitem{}
Withers, C. S. and Nadarajah, S. (2008)  The $n$th power of a matrix 
and approximations for large $n$.  Preprint http://arxiv.org/abs/0802.0502  
%{\bf 138corb }

\end{thebibliography}
\end{document}